\documentclass[10pt]{amsart}

\usepackage{amsmath,amsfonts, latexsym,graphicx, footmisc, amssymb, mathtools, enumerate}

\usepackage{hyperref}
\hypersetup{colorlinks=true, urlcolor=blue, citecolor=blue, linkcolor=blue}

\newtheorem{defn0}{Definition}[section]
\newtheorem{prop0}[defn0]{Proposition}
\newtheorem{conj0}[defn0]{Conjecture}
\newtheorem{thm0}[defn0]{Theorem}
\newtheorem{lem0}[defn0]{Lemma}
\newtheorem{corollary0}[defn0]{Corollary}
\newtheorem{example0}[defn0]{Example}
\newtheorem{remark0}[defn0]{Remark}
\newtheorem{question0}[defn0]{Question}
\newtheorem{exercise0}[defn0]{Exercise}

\title[Addendum on Aligned Hypersurface Singularities]{A Non-isolated Addendum to the Theorem \\of Bobadilla and Pe\l ka}

\author{David B. Massey}

\date{}

\begin{document}


\maketitle




In \cite{bobpelka}, Bobadilla and Pe\l ka prove the long-standing parameterized version of the Zariski Multiplicity Conjecture (ZMC) (see \cite{zarques}, \cite{eyralsurvey}): a (continuously) parameterized family of isolated hypersurface singularities with constant local, ambient topological-type has constant multiplicity. They prove this by proving that a parameterized family of isolated hypersurface singularities with constant Milnor number has constant multiplicity, and then using the well-known result of Teissier \cite{teissierdeform} (see also L\^e, \cite{leattach} and \cite{topsing}) that the Milnor number is an invariant of the local, ambient topological-type.

In \cite{levar1}, \cite{levar2}, \cite{lecycles}, we developed the L\^e numbers as a generalization of the Milnor number for hypersurfaces with non-isolated singularities. By combining the Uniform L\^e-Iomdin Formulas (Theorem 4.15 of \cite{lecycles}) with the results of Bobadilla and Pe\l ka, we easily conclude that an analytically parameterized family of hypersurfaces with constant L\^e numbers have constant multiplicity. However, the L\^e numbers are {\bf not} invariants of the local, ambient topological-type, and so we cannot conclude a non-isolated version of the ZMC in this manner.

But, in Definition 7.1 of \cite{lecycles}, we defined a class of hypersurface singularities for which the L\^e numbers are, in fact, topological invariants: {\it aligned singularities}. These are hypersurfaces that have a good ($a_f$, without Whitney conditions) stratification in which the strata contained in the singular set all have smooth closures, e.g., a hyperplane arrangement, a hypersurface with a smooth curve of singular points, or a hypersurface with a smooth 2-dimensional singular set such that the transverse Milnor number at points in the singular set jumps up only along a smooth curve. The invariance of the L\^e numbers for aligned singularities is contained in Corollaries 7.7 and 7.8 of \cite{lecycles}.

\medskip

Thus, from our proof of Theorem 7.9 in \cite{lecycles},  we are able to immediately conclude the following corollary to the theorem of Bobadilla and Pe\l ka:

\bigskip

{\bf \noindent {Corollary}}: A complex analytic parameterized family of aligned hypersurface singularities with constant local, ambient topological-type has constant multiplicity.

\bigskip

{\bf \noindent {Remark}}: The Theorem of Bobadilla and Pe\l ka is for continuous families of power series. It is unclear to us if our corollary above remains true with just the assumption of continuity in the parameter.

\bibliographystyle{plain}

\bibliography{Masseybib}

\end{document}